# Labeled Mereological Set Theory

Zuhair Al-Johar

August 5, 2012

**Introduction:** The methodology used here might provide a neat method of examining paradoxes and ways to circumvent them. Most of the known set theoretic paradoxes (Russell's, Cantor's, Burali-Forti's,..) can be paralleled here and examined. This account will shed the light on a particular application of this method that appears to elude paradoxes; an application that have shortcomings that will be illustrated here and suggestions to solutions are proposed. First I'll present the exposition of the theory, and then speak about its background and the aims behind it and how it can extend our knowledge of overcoming paradoxes.

**EXPOSITION**: Labeled Mereological Set Theory "LMST" is a theory in first order logic with primitives of identity "=", membership "∈" and label "L". L is a one place function symbol; Lx=y is read as y is the label of x.

A formula Φ is said to be *cyclic* iff a cyclic undirected multi-graph can be defined on Φ such that each node is a variable in Φ and each edge between variables a,b is an occurrence of an atomic subformula of Φ whose sole arguments are a,b.

A formula Φ is said to be *multi-cyclic* iff there are at least two occurrences of formulas in Φ each being cyclic.

*Axioms:* Identity axioms +

(1) Sets with the same members are identical.

(2) All sets have members. *

(3) A singleton is its member.

(4) A member is singleton.

(5) If Φ holds for a singleton then a set of all singletons satisfying Φ exist.

(6) Distinct sets have distinct labels.

(7) If Φ,π are parameter free non multi-cyclic formulas then:
{x| ∃y. Φ(y) ∧ x ∈ Ly} = {x| ∃y. π (y) ∧ x ∈ Ly} → (π ↔ Φ).

(8) There are distinct sets.

* The exclusion of the empty set by axiom (2) was only adopted with respect to strict adherence to the mereological background of this theory, we can actually drop it and remove the restrictive antecedent in comprehension (5) and we'll round of this method of representing formulae by object extensions by representing the never fulfilled formulas with the empty set.

I chose this method because it appears simpler to examine paradoxes with set theory through it. The set structure is very simple and actually copies basic mereological principles, principles that are well understood at informal level, Labels are actually also too simply axiomatized.

The aim is to find object extensions of predicates in such a manner that when those predicates are not equivalent then their extensions would be distinct and vise-verse.

The plan to do that is to define the set of all elements of labels of sets satisfying Φ, in symbols {x| ∃y. Φ(y) ∧ x ∈ Ly}, and that set will stand as the object extending Φ.

Now to see how this can be done in the system above: review the first six axioms above, and try to obtain a paradox in a Russell like manner. It will be seen that a paradox cannot be obtained in this way! For an example take the set X of all elements of "Labels" that are not "subsets" of what they label. Now either X has a Label that is not a subset of it but by then X will seize to be the set of elements of all Labels that are not subsets of what they label because its label is not a subset of it! so the label of X must be a subset of X which gives the apparent contradiction with the condition defining X. The solution is that the label of X is indeed a subset of X that completely overlaps with subsets of X that are labeling sets that don't have their labels as subsets of them, and there is nothing in this methodology against having such overlaps between labels, thus no paradox is raised!!! In a similar (if not much easier) manner one can purge away the apparent paradox involved with X being the set of all elements of Labels that are not "elements" of what they label.

Of course the above are just examples so we cannot make a generalization here, yet what is shown is expected paradoxical formulas not managing to evoke a paradox, similarly Burali-Forti's and Cantor's paradox are purged so is Leśniewski's paradox of singletons. So this illustrates the potential of this method in avoiding known set theoretic paradoxes.

However this method also has its shortcomings like for example we'd be losing our goal which is finding extensions of predicates the distinctiveness of which reflect non equivalence of the predicates they extend. To show this:

What we want to do is to have the following: For predicates Φ, π

{x| ∃y. Φ(y) ∧ x ∈ Ly} = {x| ∃y. π(y) ∧ x ∈ Ly} → (π ↔ Φ)

This fails for the general case of predicates Φ and π.

Proof:

Let X={y| ∃z. (¬Lz ⊂ z) ∧ y ∈ Lz}

Let X*={y| ∃z. ((¬Lz ⊂ z) ∨ (z=X)) ∧ y ∈ Lz}

Now we have X=X*! But the predicate (¬Lz ⊂ z) is not equivalent to the predicate ((¬Lz ⊂ z) ∨ (z=X)) since LX ⊂ X !

So two non equivalent predicates are extended by the SAME set.

Thus the extensional paralleling of predicates as depicted above fails!

So to make matters clear the axioms up to (6) do not cause any paradox, but it causes some "*confusion*" of extending predicates for the general case of predicates which is a shortcoming not a paradox!!!

However that the goal of this method failed for the general case doesn't mean we cannot have some particular cases where it holds, and those particular cases are worth contemplating since they may increase our insight about extending predicates.

We need to see first which sort of formulas is the cause of the confusion. Of course clearly cyclical formulas are the culprit, for it is those kinds of expressions that cause double representation as above and it is indeed the source of confusion. So if we only use Acyclic formulas then the above extensional paralleling concept would hold for all of them. However this is too harsh a measure, the culprit seems to be something less than that, actually how I see matters is if we allow in addition to acyclic formulas, other formulas that are indeed cyclic like the formula ¬Lz ⊂ z but forbid formulas that use two or more occurrences of cyclic formulas in them, then we can also resolve this situation. Of course we must forbid having parameters since this will re-introduce the confusion through parameters. So we can use parameter free cyclic formulas provided that they don't have more than one cyclic formula occurring in them. I think this resolves the matter. Accordingly I see that multi-cyclic formulas are the culprit.

Another possible line is to be less cautious and allow using even multi-cyclic parameter free formulas as long as they don't have two separate cyclic subformulas equivalent to each other; so for a formula Φ if φ1,..., φn are "all" separate occurrences of cyclic subformulas of Φ, then if ¬(φj ↔ φi) for ¬i=j, then Φ is allowed.

In this way we might have a better chance of representing more predicates by object extensions that reflect their equivalence, therefore adding more strength to comprehension over sets!

It's worth noticing that the known stratification criterion of Quine has been established to be reducible to "Acyclicity" criterion [1,2]. So the theory New Foundation "NF" and its fragments all depends on the strict notion of acyclicity, while this method goes beyond that to encounter frank violations of acyclicity (and stratification) in parameter "free" formulas yet without apparently encountering any paradox, this will serve to widen our definitional coverage over sets i.e. comprehension. So this method is one step ahead of acyclicity and stratification and it might increase our understanding of sets.

**References:**

[1] Al-Johar Z: Acyclic Comprehension is equal to Stratified Comprehension. http://zaljohar.tripod.com/acycliccomp.pdf

[2] Al-Johar Z; Holmes M.R; Bowler N: The Axiom Scheme of Acyclic Comprehension. http://math.boisestate.edu/~holmes/acyclic_abstract_final.pdf